\documentclass[12pt]{article}
\usepackage{graphicx,amsmath,amsthm,amssymb,dsfont, mathrsfs }
\numberwithin{equation}{section}
\def\NN{\mathbb N}

\def\ZZ{\mathbb Z}

\def\FF{\mathbb F}

\def\mk{\boldsymbol{\mathsf{k}}}
\def\mr{\boldsymbol{\mathsf{r}}}

\def\msu{\boldsymbol{\mathsf{u}}}
\def\sa{\mathsf a}

\def\fT{\mathfrak T}

\def\bL{{\bf L}}

\def\bT{{\bf T}}
\def\bU{{\bf U}}

\def\bx{{\bf x}}
\def\by{{\bf y}}

\def\ba{{\bf a}}
\def\bl{{\bf l}}

\def\bk{{\bf k}}
\def\bv{{\bf v}}

\def\bu{{\bf u}}
\def\bw{{\bf w}}

\def\bgamma{\boldsymbol{\gamma}}

\def\bs{{\bf 0}}
\def\b1{{\bf 1}}
\def\d1{\mathds{ 1}}

\def\mod{{\rm mod}}
\def\ad{{\rm and}}
\def\oor{{\rm or}}
\def\card{{\rm card}}
\def\with{{\rm with}}
\def\where{{\rm where}}
\def\for{{\rm for}}
\def\fall{{\rm for\; all}}
\def\vol{{\rm vol}}
\def\Re{{\rm Re}}

\def\Tr{{\rm Tr}}

\def\qed{ \ \vrule width.2cm height.2cm depth0cm\smallskip}

\begin{document}



\title{On a bounded remainder set for a digital Kronecker sequence}
\author{Mordechay B. Levin}

\date{}

\maketitle

\begin{abstract}
Let $\bx_0,\bx_1,...$ be a sequence of points in $[0,1)^s$.
 A subset $S$ of $[0,1)^s$ is called a bounded remainder set if there exist two real numbers $a$ and $C$ such that, for every integer $N$,
$$
  | {\rm card}\{n <N \; | \;   \bx_{n} \in S  \} -  a N| <C  .
$$

Let $ (\bx_n)_{n \geq 0} $ be an $s-$dimensional digital Kronecker-sequence in base $b \geq 2$,
$\bgamma =(\gamma_1,...,\gamma_s)$,
$\gamma_i \in [0, 1)$ with $b$-adic expansion\\ $\gamma_i= \gamma_{i,1}b^{-1}+ \gamma_{i,2}b^{-2}+...$, $i=1,...,s$.
 In this paper, we prove that $[0,\gamma_1) \times ...\times [0,\gamma_s)$ is the bounded remainder set  with respect to the
sequence $(\bx_n)_{n \geq 0}$ if and only if
\begin{equation} \nonumber
   \max_{1 \leq i \leq s} \max \{ j \geq 1 \; | \; \gamma_{i,j} \neq 0 \}  < \infty.
\end{equation}
\end{abstract}
%
%
\section{Introduction }
{\bf 1.1. Discrepancy.} Let $\bx_0,\bx_1,...$ be a sequence of points in $[0,1)^s$,  $S \subseteq [0,1)^s$,
\begin{equation}\label{In1}
\Delta(S, (\bx_{n})_{n=0}^{N-1}  )= \sum_{n=0}^{N-1}  ( \b1_{S}(\bx_{n}) - \lambda(S)),
\end{equation}
where $\d1_{S}(\bx) =1, \; {\rm if} \;\bx  \in S$,
and $   \d1_{S}(\bx) =0,$  if $
\bx   \notin S$.  Here $\lambda(S)$ denotes the $s$-dimensional Lebesgue-measure of $S$.
We define the star {\it discrepancy} of an
$N$-point set $(\bx_{n})_{n=0}^{N-1}$ as
\begin{equation} \label{In2}
   \emph{D}^{*}((\bx_{n})_{n=0}^{N-1}) =
    \sup\nolimits_{ 0<y_1, \ldots , y_s \leq 1} \; |
  \Delta([\bs,\by),(\bx_{n})_{n=0}^{N-1})/N |,
\end{equation}
where $[\bs,\by)=[0,y_1) \times \cdots \times [0,y_s) $.
 The sequence $(\bf x_n)_{n \geq 0}$ is said to be
{\it uniformly distributed}  in $[0,1)^s$
if $D_N \to \; 0$.

In 1954, Roth proved that
$
   \limsup\nolimits_{N \to \infty } N (\ln N)^{-\frac{s}{ 2}} \emph{D}^{*}((\bx_{n})_{n=0}^{N-1})>0 . $
According to the well-known conjecture (see, e.g., [BeCh, p.283]), this estimate can be improved
 to
$
 \limsup\nolimits_{N \to \infty } N (\ln N)^{-s} \emph{D}^{*}((\bx_{n})_{n=0}^{N-1})>0
$.
See [Bi] and [Le1] for the  results on this conjecture.

 An $s$-dimensional sequence  $((\bx_{n})_{n \geq 0})$ is of
 \texttt{low discrepancy} (abbreviated
l.d.s.) if $ \emph{D}^{*}((\bx_{n})_{n=0}^{N-1})=O(N^{-1}(\ln
N)^{s}) $ for $ N \rightarrow \infty $.
For  examples of  l.d.s. see, e.g.,  in [BeCh], [DiPi], [Ni]. \\ \\
{\bf 1.2. Digital Kronecker sequence.}

For an arbitrary prime power $b$, let $\FF_b$ be the finite field of order $b$,  $\FF_b^{*} =\FF_b\setminus 0$, $\ZZ_b = \{0, 1, . . . , b-1 \}$. Let $\FF_b[z]$ be the
set of all polynomials over $\FF_b$,
and let $\FF_b((z^{-1}))$ be the field of formal Laurent series.
Every element $L$
of $\FF_b((z^{-1}))$ has a unique expansion into a formal Laurent series
\begin{equation}   \label{In10}
   L = \sum_{k=w}^{\infty} u_k z^{-k}  \quad \with \quad u_k \in \ZZ_b, \;\; \and \;\; w\in \ZZ
	  \;\; \where \;\; u_w  \neq 0.
\end{equation}
The discrete exponential evaluation $\nu$ of $L$ is defined by
\begin{equation} \nonumber
\nu(L) := −w, \quad L\neq 0, \qquad \nu(0) := −\infty.
\end{equation}
Furthermore, we define the ``fractional part`` of $L$ by
\begin{equation}   \label{In14}
   \{L\} = \sum_{k= \max(1,w)}^{\infty} u_k z^{-k} .
\end{equation}
We choose bijections  $\psi_r : \ZZ_b \to  \FF_b $ with  $\psi_r(0) = 0$, and for $i = 1, 2, ... , s$
and $j = 1, 2, . . . $ we choose bijections $\eta_{i,j} : \FF_b \to  \ZZ_b $.
For n = 0, 1, ..., let
\begin{equation}   \label{In16}
  n =\sum_{r=0}^{\infty} a_r(n) b^r
\end{equation}
be the digit expansion of $n$ in base $b$, where $a_r(n) \in \ZZ_b$ for $r \geq 0$ and
$a_r(n) = 0$ for all sufficiently large $r$.

 With every n = 0, 1, . . ., we associate the polynomial
\begin{equation}   \label{In18}
  n(z) =\sum_{r=0}^{\infty} \psi_r(a_r(n)) z^r \in \FF_b[z]
\end{equation}
and if $L  \in \FF_b((z^{-1}))$ is as in (\ref{In10}), then we define
\begin{equation}   \label{In20}
   \eta^{(i)}(L) = \sum_{k= \max(1,w)}^{\infty} \eta_{i,k}(u_k) b^{-k} .
\end{equation}

In [Ni], Niederreiter introduced a non-Archimedean analogue of the classical Kronecker sequences.
For every $s$-tuple $\bL = (L_1, . . . , L_s)$ of elements of $\FF_b((z^{-1}))$, we define the sequence
$S(\bL) = (\bl_n)_{n \geq 0}$ by
\begin{equation}   \label{In22}
\bl_n = (l^{(1)}_n,..., l^{(s)}_n   ), \quad
l^{(i)}_n = \eta^{(i)}( n(z) L_i(z)),
 \quad \for \quad  \; 1 \leq i \leq s, \; n \geq 0.
\end{equation}
This sequence is sometimes called a digital Kronecker sequence (see [LaPi, p.4]).
In analogy to classical Kronecker sequences, in [LaNi, Theorem 1], the following theorem  has been
 proven \\ \\
 {\bf  Theorem A.} {\it A digital
Kronecker sequence $S(L)$ is uniformly distributed in $[0,1)^s$ if and only if $1, L_1, . . . ,
L_s$ are linearly independent over $\ZZ_b[x]$.}\\


 By $\mu_1$ we denote the normalized Haar-measure
on $\FF_b((z^{-1}))$ and by $\mu_s$ the $s$-fold product measure on $\FF_b((z^{-1}))^s$.
In [La1], Larcher proved the following metrical upper bound on the star discrepancy
of digital Kronecker sequences $ D_N(S(L)) = O(
N^{-1}(\log N)^s(\log \log N)^{2+\epsilon})$.
For $\mu_s$-almost all $L \in \FF_b((z^{-1}))^s$, $\epsilon >0$.

In  [LaPi, p.4], Larcher and Pillichshammer were able to give corresponding metrical
lower bounds for the discrepancy of digital Kronecker
sequences $ D_N(S(L))  \geq
c(b, s) N^{-1}
(\log N)^s \log \log N $
 for $\mu_s$-almost all $L \in \FF_b((z^{-1}))^s$, for infinitely many $N \geq 1$
with some $c(b, s) > 0$ not depending on N. \\ \\
{\bf 1.3. Bounded remainder set}.\\
 {\bf Definition 1}. {\it Let $\bx_0,\bx_1,...$ be a sequence of point in $[0,1)^s$.
 A subset $S$ of $ [0,1)^s$ is called a {\sf bounded remainder set} for $(\bx_n)_{n \geq 0}$
if the discrepancy function $\Delta(S, (\bx_n)_{n = 0}^{N-1})$ is bounded in $N$.}\\


 Let $\alpha$ be an irrational number, let I be an interval in $[0,1)$ with the length $|I|$, let $\{n\alpha\}$ be the fractional part of $n\alpha$, $n=1,2,…$.
 Hecke,  Ostrowski  and  Kesten  proved that $\Delta(S, (\{n\alpha\})_{n = 1}^N)$ is bounded
 if and only if $|I|=\{k\alpha\}$ for some integer $k$ (see references in [GrLe]).

 The sets of bounded remainder for the classical $s$-dimensional Kronecker sequence were studied
by Lev and Grepstad [GrLe]. The case of Halton's sequence was studied by Hellekalek [He].
For references to others investigations on   bounded remainder set
 see [GrLe].

Let $\bgamma =(\gamma_1,...,\gamma_s)$,
$\gamma_i \in (0, 1)$ with $b$-adic expansion  $\gamma_i= \gamma_{i,1}b^{-1}+ \gamma_{i,2}b^{-2}\\
+...$, $i=1,...,s$.
 In this paper, we prove\\ \\
{\bf Theorem.} {\it Let $(\bl_n)_{n \geq 0}$ be a uniformly distributed digital Kronecker sequence.
The set $[0,\gamma_1) \times ...\times [0,\gamma_s)$ is of bounded remainder  with respect to
  $(\bl_n)_{n \geq 0}$ if and only if }
\begin{equation} \label{Cond}
   \max_{1 \leq i \leq s} \max \{ j \geq 1 \; | \; \gamma_{i,j} \neq 0 \}  < \infty.
\end{equation}

In [Le2], we proved similar results for digital sequences described in [DiPi, Sec. 8]. Note that according to
Larcher's conjecture [La2, p.215], the assertion of the Theorem is true for all digital $(t,s)$-sequences in
base $b$.

\section{Notations.}


A subinterval $E$ of $[0,1)^s$  of the form
$$  E = \prod_{i=1}^s [a_ib^{-d_i},(a_i+1)b^{-d_i}),   $$
 with $a_i,d_i \in \ZZ, \; d_i \ge 0, \; 0  \le a_i < b^{d_i}$ for $1 \le i \le s$ is called an
{\it elementary interval in base $b \geq 2$}.\\ \\
{\bf Definition 2}. {\it Let $0 \le t \le m$  be  integers. A {\sf $(t,m,s)$-{\sf net in base $b$}} is a point set
$\bx_0,...,\bx_{b^m-1}$ in $ [0,1)^s $  such that $\# \{ n \in [0,b^m -1] | x_n \in E \}=b^t$   for every elementary interval E in base  $b$ with
$\vol(E)=b^{t-m}$.}\\  \\
{\bf Definition 3.} (\cite[Definition 4.30]{DiPi}) { \it For a given dimension $s \geq 0$, an integer base $b \geq 2$, and a
function $\bT : \NN_0 \to \NN_0$ with $\bT(m) \leq  m$ for all $m \in \NN_0$, a sequence $(\bx_0,\bx_1, . . .)$
of points in $[0, 1)^s$ is called a $(\bT, s)$-sequence in base $b$ if for all integers $m \geq 1$
and $k \geq  0$, the point set consisting of the points $x_{kb^m}, . . . ,x_{kb^m+b^m-1}$ forms
a $(\bT(m),m, s)$-net in base $b$. }\\

 A $(\bT, s)$-sequence in base $b$ is called a strict $(\bT, s)$-sequence in
base $b$ if for all functions $\bU : \NN_0 \to \NN_0$ with $\bU(m) \leq m$ for all $m \in \NN_0$ and with
$\bU(m) < \bT(m)$ for at least one $m \in \NN_0$, it is not a $(\bU, s)$-sequence in base $b$.\\ \\
{\bf Definition 4.} ([DiNi, Definition 1]) { \it
Let $m, s \geq 1$ be integers. Let $C^{(1,m)},...,$ $C^{(s,m)}$ be $m \times m$ matrices over $\FF_b$.
Now we construct $b^m$ points in $[0, 1)^s$.
 For $ n= 0, 1,...,b^m-1$, let $n =\sum^{m-1}_{j=0} a_j(n) b^{j}$
be the $b$-adic expansion of $n$.  For $r = 0,1,...$
we choose bijections  $\psi_r : \ZZ_b \to  \FF_b $ with  $\psi_r(0) = 0$, and for $i = 1, 2, . . . , s$
and $j = 1, 2, . . . $ we choose bijections $\eta_{i,j} : \FF_b \to  \ZZ_b $.
We map the vectors
\begin{equation} \label{No1}
	y^{(i,m)}_{n}=(y^{(i,m)}_{n,1},...,y^{(i,m)}_{n,m}),\quad y^{(i,m)}_{n,j}=\sum_{r=0}^{m-1} \psi_r(a_r(n)) c^{(i,m)}_{j,r}\in  \FF_b
\end{equation}
to the real numbers
\begin{equation}   \label{No2}
   x^{(i)}_n =\sum_{j=1}^m \eta_{i,j} (y^{(i,m)}_{n,j})/b^j
\end{equation}
to obtain the point
\begin{equation} \nonumber
   \bx_n= (x^{(1)}_n,...,x^{(s)}_n) \in [0,1)^s.
\end{equation}
 \\

The point set  $ \{\bx_0,...,\bx_{b^m-1} \}$ is called a {\sf digital net} (over $\FF_b$) (with {\sf generating matrices} $(C^{(1,m)},...,C^{(s,m)}) $).

For $m = \infty$, we obtain a sequence $\bx_0, \bx_1,...$ of   points in $[0, 1)^s$  which is called a {\sf digital sequence} $($over $\FF_b)$ $($with {\sf generating matrices} $(C^{(1,\infty)},...,C^{(s,\infty)}) )$.}

We abbreviate  $C^{(i,m)}$ as $C^{(i)}$ for $m \in \NN$ and for $m=\infty$. \\ \\
 {\bf Lemma A} ([LaNi, ref. 1-8]).     { \it A digital Kronecker sequence in base $b$ can be expressed as some  digital $(\bT, s)$-sequence in base $b$.}\\ \\
 {\bf Lemma B} (\cite[Theorem 4.86]{DiPi}). {\it Let $b$ be a prime power. A strict digital
$(\bT, s)$-sequence over $\FF_b$ is
uniformly distributed modulo one, if and only if $\liminf_{m \to \infty} (m - \bT(m))=\infty$.}\\

For $m>n$, we put $\sum_{j=m}^n c_j= 0$ and  $\prod_{j=m}^n c_j =1$.
For $x =\sum_{j \geq 1}  x_{j} b^{-j}$,
where $x_{i} \in  \ZZ_b =\{0,...,b-1\}$,   we define the truncation
\begin{equation}  \nonumber
        [x]_m =\sum_{1 \leq j \leq m}  x_{j}b^{-j} \quad \with \quad m \geq 1.
\end{equation}
If $\bx = (x^{(1)}, . . . , x^{(s)})  \in [0, 1)^s$, then the truncation $[\bx]_m$ is defined coordinatewise, that is, $[\bx]_m = ( [x^{(1)}]_m, . . . , [x^{(s)}]_m)$.

 For $x =\sum_{j \geq 1}  x_{j} b^{-j}$ and $y =\sum_{j \geq 1}  y_{j}b^{-j}$
where $x_j ,y_j \in  \ZZ_b$, we define the ($b$-adic) digital shifted point $v$ by
$v = x \oplus y := \sum_{j \geq 1}  v_{j}b^{-j}$,
 where $v_j \equiv x_j + y_j \;(\mod \;b)$ and $v_j \in \ZZ_b$.
For $\bx =(x^{(1)}, . . . , x^{(s)}) \in [0, 1)^s$ and $\by = (y^{(1)}, . . . , y^{(s)}) \in [0, 1)^s$, we define the ($b$-adic) digital shifted point $\bv$ by
$ \bv =\bx \oplus \by =(x^{(1)} \oplus y^{(1)}, . . . ,x^{(s)} \oplus y^{(s)} )$.
 For $n_1,n_2 \in [0,b^m)$, we define
$n_1 \oplus n_2 := (n_1 /b^m\oplus n_2)b^m)b^m$.

For $x =\sum_{j \geq 1}  x_{j} b^{-j}$,
where $x_{j} \in  \ZZ_b$,  $x_j=0 $ $(j=1,...,k)$ and $x_{k+1} \neq 0$, we define the
absolute valuation  $\left\|.  \right\|_b $ of $x$ by  $\left\|x  \right\|_b =b^{-k-1}$.
Let $\left\| n  \right\|_b =b^k$ for $n \in [b^k,b^{k+1})$.\\ \\
{\bf Definition 5.} {\it A sequence   $ (\bx_{n})_{n \geq 0} $
	in $[0,1)^s$ is  {\sf weakly admissible} in
base $b$  if}
\begin{equation}  \label{3}
  \varkappa_m:= \min_{0 \leq k <n < b^m} \left\| \bx_n \ominus \bx_k  \right\|_b
  >  0\quad \forall m \geq 1\; {\rm where} \;\; \left\| \bx  \right\|_b := \prod_{i=1}^s
	\left\|x^{(i)}  \right\|_b .
\end{equation}

Let $p$ be a prime, $b=p^{\kappa}$,
\begin{equation} \nonumber
 E(\alpha) := exp (2\pi i \Tr(\alpha)/p), \qquad \alpha \in \FF_b,
\end{equation}
where $\Tr : \FF_b \to  \FF_p$ denotes the usual trace of an element of $\FF_b$ in $\FF_p$.

Let
\begin{equation} \label{Del1}
   \delta(\fT) =   \begin{cases}
    1,  & \; {\rm if}  \;  \fT  \;{\rm is \;true},\\
    0, &{\rm otherwise}.
  \end{cases}
\end{equation}\\

By [LiNi, ref. 5.6 and  ref. 5.8], we get
\begin{equation}   \label{Del2}
 \frac{1}{q} \sum_{\beta \in \FF_b} E(\alpha \beta) = \delta (\alpha =0 ), \quad \where \quad
  \alpha  \in \FF_b.
\end{equation}

\section{Proof}
{\bf Lemma 1.} {\it Let  $(\bx_n)_{n \geq 0}$ be a  weakly admissible digital sequence in base $b$,
  $m \geq 1$, $\tau_m =[\log_b (\kappa_m)]+m $. Then we have for all integers $A \geq 1$}
\begin{equation} \nonumber
    | \Delta([\bs,\bgamma),(\bx_{n})_{n=b^m A}^{b^m A+N-1}) -
  \Delta([\bs,[\bgamma]_{\tau_m}),(\bx_{n})_{n=b^m A}^{b^m A+N-1})| \leq s,\quad
	  \forall N\in[1,b^m].
\end{equation}\\
{\bf Proof.} Let
\begin{equation} \nonumber
   B=[\bs,\bgamma), \quad  B_i =\prod_{1 \leq j \leq s, j\neq i} [0, \gamma^{(j)})
	\times [0, [\gamma^{(i)}]_{\tau_m}) \quad \ad \quad B_0 =\cup_{i=1}^s(B \setminus B_i).
\end{equation}
It is easy to see that $B= [\bs,[\bgamma]_{\tau_m}) \cup B_0$.
 By (\ref{In1}), we get
\begin{equation} \nonumber
  \Delta([\bs,\bgamma),(\bx_{n})_{n=b^m A}^{b^m A+N-1})
\end{equation}
\begin{equation} \nonumber
      = \Delta([\bs,[\bgamma]_{\tau_m}),(\bx_{n})_{n=b^m A}^{b^m A+N-1})  +\Delta(B_0,(\bx_{n})_{n=b^m A}^{b^m A+N-1}) .
\end{equation}
Hence
\begin{equation} \nonumber
     |\Delta([\bs,\bgamma),(\bx_{n})_{n=b^m A}^{b^m A+N-1}) -
  \Delta([\bs,[\bgamma]_{\tau_m}),(\bx_{n})_{n=b^m A}^{b^m A+N-1}) |
\end{equation}
\begin{equation} \label{Prof1}
  \leq \sum_{i=1}^s  | \Delta(B \setminus B_i,(\bx_{n})_{n=b^m A}^{b^m A+N-1})| .
\end{equation}
Suppose that there exist $i \in [1,s]$,  $ k,n \in [0, b^m),\; k \neq n$ and $A \geq 1$
such that  $x_{n+ b^m A}, x_{k+ b^m A} \in B \setminus B_i$.
Therefore
\begin{equation} \nonumber
     x_{n+ b^m A,j}^{(i)} = x_{k+ b^m A,j}^{(i)} \quad \for \quad j=1,...,\tau_m.
\end{equation}
From (\ref{In16}),  (\ref{No1}) and \eqref{No2},   we have
\begin{equation} \nonumber
     y_{n+ b^m A,j}^{(i)} = y_{k+ b^m A,j}^{(i)} \quad \for \quad j=1,...,\tau_m,
\end{equation}
\begin{equation} \nonumber
     y_{n+ b^m A,j}^{(i)}=y_{n,j}^{(i)}+y_{b^m A,j}^{(i)}, \quad \ad \quad
		 y_{k+ b^m A,j}^{(i)}=y_{k,j}^{(i)}+y_{b^m A,j}^{(i)} \quad \for \; j=1,...,\tau_m.
\end{equation}
Hence
\begin{equation} \nonumber
     y_{n,j}^{(i)} = y_{k,j}^{(i)}, \;\;  j=1,...,\tau_m  \quad \ad \quad
     x_{n,j}^{(i)} = x_{k,j}^{(i)}, \;\;  j=1,...,\tau_m.
\end{equation}
Therefore
\begin{equation} \nonumber
   \left\| x_{n}^{(i)} \ominus x_{k }^{(i)}\right\|_b < b^{-\tau_m}\leq \kappa_m \quad \ad \quad
   \left\| \bx_{n} \ominus \bx_{k}\right\|_b \geq \varkappa_m.
\end{equation}
By (\ref{3}) we have a contradiction. Thus
\begin{equation} \nonumber
\card\{n \in [0,b^m) \;|\; \bx_{n+ b^m A} \in B \setminus B_i \} \leq 1,\; \ad \;
 | \Delta(B \setminus B_i,(\bx_{n})_{n=b^m A}^{b^m A+N-1})| \leq 1.
\end{equation}
Using (\ref{Prof1}),  we get the assertion of Lemma 1. \qed  \\


 Let $\beta_1,...,\beta_{\kappa}$ be a $F_p$ basis of $\FF_b$, and let $\Tr$ be a standard trace function.
	Let
	\begin{equation}  \label{Lemm1}
  \omega(\alpha) =\sum_{j=1}^{\kappa} p^{j-1} \Tr(\alpha \beta_j), \qquad b= p^{ \kappa}.
\end{equation}	
We use notations  (\ref{In16}), (\ref{No1}) and (\ref{No2}).
Let $ n =\sum_{r\geq 0} a_r(n) b^{r}$	be the $b$-adic expansion of $n$, and let
\begin{equation}  \label{Lemm2}
   \tilde{n} =\sum_{r\geq 0} \omega(\psi_r(a_r(n))) b^{r}.
\end{equation}
Therefore
\begin{equation} \label{Lemm26a}
\{ \tilde{n}\;| \;0 \leq n <b^m\} = \{0,1,...,b^m-1 \}.
\end{equation}	
Hence	
\begin{equation}  \label{Lemm2a}
  \psi_r(a_r(n)) =   \omega^{-1}(a_r(\tilde{n}))
\end{equation}	
and
\begin{equation}  \label{Lemm3}
    u^{(i)}_{\tilde{n},j}:=
	            \sum_{r \geq 0}  \omega^{-1}(a_r(\tilde{n})) c^{(i)}_{j,r}
	=  \sum_{r \geq 0} \psi_r(a_r(n)) c^{(i)}_{j,r}  = y^{(i)}_{n,j}
								 ,\quad 1 \leq i \leq s.
\end{equation}
Let
\begin{equation} \label{Lemm4}
x^{(s+1)}_{n} =\{ n/b^m \}, \;\; x^{(s+1)}_{n,j}=a_{m-j}(n),\;\;y^{(s+1)}_{n,j}=\psi_{m-j}(x^{(s+1)}_{n,j}).
\end{equation}
Bearing in mind that $  \psi_{m-j}(a_{m-j}(n)) =   \omega^{-1}(a_{m-j}(\tilde{n})) $, we put
\begin{equation} \label{Lemm5}
	u^{(s+1)}_{\tilde{n},j}:=    \omega^{-1}(a_{m-j}(\tilde{n}))=\psi_{m-j}(a_{m-j}(n))=y^{(s+1)}_{n,j},
 \;\;                       j\in[1,m].
\end{equation}
Let
%
\begin{equation} \nonumber
	      u^{(i)}_{n} =(u^{(i)}_{n,1},...,u^{(i)}_{n,\tau_m}) \in \FF_b^{\tau_m} \quad \ad \quad
           u^{(s+1)}_{n} =(u^{(s+1)}_{n,1},...,u^{(s+1)}_{n,m}).
\end{equation}
We abbreviate  $s+1$-dimensional vectors
$(u^{(1)}_{n},...,u^{(s+1)}_{n})$, $(k^{(1)},..., k^{(s+1)})$ and $(r^{(1)},..., r^{(s+1)})$
 by symbols  $\msu_{n}$, $\mk$ and $\mr$, and $s$-dimensional vectors
$(u^{(1)}_{n},...,\\u^{(s)}_{n})$, $(k^{(1)},..., k^{(s)})$
 by symbols  $\bu_{n}$ and $\bk$. \\

By  (\ref{Lemm1}) - (\ref{Lemm5}), we get  $u^{(s+1)}_{n} =  u^{(s+1)}_{n +b^m A}     $, $A=1,2, ...$,
\begin{equation} \label{Lemm6}
	u^{(i)}_{n_1 \oplus n_2,j} = u^{(i)}_{n_1,j}+ u^{(i)}_{n_2,j},\; j \geq 1, i \in[1, s+1], \;\;	\msu_{n_1 \oplus n_2} = \msu_{n_1}+ \msu_{n_2}.
\end{equation}


 Let  $N \in [1,b^m]$, $\gamma^{(s+1)}=N /b^m $,	$k  = \sum_{j = 1}^{\tau_m} k_{j} b^{-j} >0$, with
$ k_j \in \ZZ_b$,
\begin{equation} \label{Lemm7}
v(k):= \max \{ j \in [1, \tau_m] \; | \; k_j \neq 0\}, \quad  v(0)=0.
\end{equation}

Similarly to [Ni, Theorem 3.10] (see also [DiPi, Lemma 14.8]), we consider the following
 Fourier series decomposition of the discrepancy function :\\ \\
%
{\bf Lemma 2.} {\it Let  $A \geq 1$ be an integer, $N \in [1,b^m]$, $\gamma^{(s+1)}=N /b^m $,
and let $(\bx_n)_{n \geq 0}$ be a digital sequence in base $b$.
Then}
\begin{equation} \nonumber
     \Delta([\bs,[\bgamma]_{\tau_m}),(\bx_{n})_{n=b^m A}^{b^m A+N-1})
\end{equation}
\begin{equation} \nonumber
 =
  \sum_{n=0}^{b^m-1} 		\sum_{(k^{(1)},...,k^{(s)})  \in \FF_b^{s\tau_m}} \; \sum_{k^{(s+1)}  \in \FF_b^{m}}
	 E(\mk \cdot \msu_{\widetilde{n+b^m A}}) \hat{\d1}(\mk)
		-b^m \prod_{i=1}^{s+1} [\gamma^{(i)}]_{\tau_m}  ,
\end{equation}
where
\begin{equation} \label{Lemm12}
 \mk \cdot \msu_n = \sum_{i=1}^{s} \sum_{j=1}^{\tau_m} k_{j}^{(i)}
u^{(i)}_{n,j} + \sum_{j=1}^m k_{j}^{(s+1)} u^{(s+1)}_{n,j},
\;\; \;\;  \hat{\d1}(\mk) = \prod_{i=0}^{s+1}   \hat{\d1}^{(i)}(k^{(i)}),
\end{equation}
 $\hat{\d1}^{(i)}(0) =[\gamma^{(1)}]_{\tau_m} \; (1 \leq i \leq s)$, $\hat{\d1}^{(s+1)}(0) =\gamma^{(s+1)}\;\;$
   and
\begin{align}  \nonumber
&   \hat{\d1}^{(i)}(k) =  b^{-v(k)}   E\Big(-\sum_{j=1}^{v(k)-1} k_j \eta_{i,j}^{-1}(\gamma_j^{(i)}) \Big)
\Big( \sum_{b=0}^{\gamma_{v(k)}^{(i)} -1}
	  E(-k_{v(k)}     \eta^{-1}_{i,v(k)}(b) )\nonumber \\
&+ E(-k_{v(k)}\eta^{-1}_{i,v(k)}(\gamma_{v(k)}^{(i)}))
		\{b^{v(k)}[\gamma^{(i)}]_{\tau_m} \}\Big), \qquad \; i \in [1,s], \nonumber \\
&    \hat{\d1}^{(s+1)}(k) = b^{-v(k)}
 E\Big( -\sum_{j=1}^{{v(k)}-1} k_j\psi_{j}(\gamma_{j}^{(s+1)})\Big)  \Big( \sum_{b=0}^{\gamma_{v(k)}^{(s+1)} -1}
	  E(-k_{v(k)}\psi_{v(k)}(b) )  \nonumber \\
&
	+ E(-k_{v(k)} \psi_{v(k)}(\gamma_{v(k)}^{(s+1)}))
		\{b^{v(k)}\gamma^{(s+1)} \} \Big). \label{Lemm12a}
\end{align}\\
{\bf Proof.}
Let $\gamma = \sum_{j = 1}^{\dot{m}} \gamma_{j} b^{-j} >0$, $w = \sum_{j = 1}^{\dot{m}} w_{j} b^{-j} $, with
$ \gamma_j,w_j \in \ZZ_b$.
It is easy to verify (see also [Ni, p. 37,38]) that
\begin{equation} \nonumber
   \d1_{[0,\gamma)}(w) =
	\sum_{r=1}^{\dot{m}} \sum_{b=0}^{\gamma_r -1} \prod_{j=1}^{r-1} \delta(w_i=\gamma_i)
	\delta(w_r=b).
\end{equation}	
By (\ref{No2}) and (\ref{Lemm3}), we have that
\begin{equation} \nonumber
   x^{(i)}_{j,n} =b   \Longleftrightarrow   y^{(i)}_{j,n} =\eta^{-1}_{i,j}(b)
	\Longleftrightarrow 	u^{(i)}_{j,\tilde{n}} =  \eta^{-1}_{i,j}(b),
\end{equation}	
and
\begin{equation} \nonumber
   \d1_{[0, [\gamma^{(i)}]_{\tau_m})}( x^{(i)}_{n}) =
	\sum_{r=1}^{\tau_m} \sum_{b=0}^{\gamma_r^{(i)} -1} \prod_{j=1}^{r-1}
	\delta(x^{(i)}_{j,n}=\gamma_j^{(i)})
	  	\delta(x^{(i)}_{r,n}=b)
\end{equation}
\begin{equation} \nonumber
  =
	\sum_{r=1}^{\tau_m} \sum_{b=0}^{\gamma_r^{(i)} -1} \prod_{j=1}^{r-1}
	\delta(y^{(i)}_{j,n}=\eta_{i,j}^{-1}(\gamma_j^{(i)}))
	  	\delta(y^{(i)}_{r,n}=\eta_{i,r}^{-1}(b))
\end{equation}	
\begin{equation} \nonumber
  =
	\sum_{r=1}^{\tau_m} \sum_{b=0}^{\gamma_r^{(i)} -1} \prod_{j=1}^{r-1}
	\delta(u^{(i)}_{j,\tilde{n}}=\eta_{i,j}^{-1}(\gamma_j^{(i)}))
	  	\delta(u^{(i)}_{r,\tilde{n}}=\eta_{i,r}^{-1}(b)), \quad i=1,...,s.
\end{equation}	
Similarly, we derive	
\begin{equation} \label{Lemm12c}
  \d1_{[0, \gamma^{(s+1)})}( x^{(s+1)}_{n}) =
	\sum_{r=1}^m \sum_{b=0}^{\gamma_r^{(s+1)} -1} \prod_{j=1}^{r-1}
	\delta(u^{(s+1)}_{j,\tilde{n}}=\psi_{j}^{-1}(\gamma_j^{(s+1)}))
	  	\delta(u^{(s+1)}_{r,\tilde{n}}=\psi_{r}^{-1}(b)) .
\end{equation}		
Let $ k \cdot u^{(i)}_{\tilde{n}} =\sum_{j=1}^{\tau_m} k_j u^{(i)}_{\tilde{n},j}$.
  By (\ref{Del2}), we have
\begin{align} \nonumber
 &  \d1_{[0, [\gamma^{(i)}]_{\tau_m})}( x^{(i)}_{n}) =
	\sum_{r=1}^{\tau_m} \sum_{b=0}^{\gamma_r^{(i)} -1} b^{-r} \sum_{k_1,...,k_r \in \FF_b} \dot{\d1}^{(i)}(k), \quad
\where   \\
& \dot{\d1}^{(i)}(k)=
	E\Big(\sum_{j=1}^{r-1} k_j (u^{(i)}_{j,\tilde{n}}-\eta_{i,j}^{-1}(\gamma_j^{(i)}))
     +	  	k_r(u^{(i)}_{r_i,\tilde{n}}-\eta_{i,r}^{-1}(b)) \Big) =  E(k \cdot   u^{(i)}_{\tilde{n}})
     \tilde{\d1}^{(i)}(k) \nonumber \\
&  \with \quad \tilde{\d1}^{(i)}(k)=   E\Big(-\sum_{j=1}^{r-1} k_j \eta_{i,j}^{-1}(\gamma_j^{(i)})
     -	  	k_r  \eta_{i,r}^{-1}(b) \Big).   \label{Lemm12aa}
\end{align}
%
%
Hence
\begin{align} \nonumber
 &  \d1_{[0, [\gamma^{(i)}]_{\tau_m})}( x^{(i)}_{n})= \sum_{r=1}^{\tau_m} \sum_{b=0}^{\gamma_r^{(i)} -1} b^{-r}
 \sum_{k_1,...,k_{\tau_m} \in \FF_b}   \delta( v(k)  \leq r)
            E(k \cdot   u^{(i)}_{\tilde{n}}) \tilde{\d1}^{(i)}(k)  \\
 &=  \sum_{k_1,...,k_{\tau_m} \in \FF_b}
	\sum_{r=1}^{\tau_m} \sum_{b=0}^{\gamma_r^{(i)} -1} b^{-r} \delta( v(k)  \leq r)
            E(k \cdot   u^{(i)}_{\tilde{n}}) \tilde{\d1}^{(i)}(k)  \nonumber  \\
& =  \sum_{k_1,...,k_{\tau_m} \in \FF_b}  E(k \cdot   u^{(i)}_{\tilde{n}}) \ddot{\d1}^{(i)}(k), \quad
 \where \quad	\ddot{\d1}^{(i)}(k) =\sum_{r= v(k)}^{\tau_m} \sum_{b=0}^{\gamma_r^{(i)} -1} b^{-r}
             \tilde{\d1}^{(i)}(k)  . \nonumber
\end{align}
Applying \eqref{Lemm12a} and \eqref{Lemm12aa}, we derive
\begin{align} \nonumber
& \ddot{\d1}^{(i)}(k) =\sum_{r=v(k)}^{\tau_m} \sum_{b=0}^{\gamma_r^{(i)} -1} b^{-r}
     E\Big(-\sum_{j=1}^{r-1} k_j \eta_{i,j}^{-1}(\gamma_j^{(i)})
     -	  	k_r  \eta_{i,r}^{-1}(b) \Big)         \\
&   =       \sum_{b=0}^{\gamma_{ v(k)}^{(i)} -1} b^{-v(k)}
     E\Big(-\sum_{j=1}^{v(k)-1} k_j \eta_{i,j}^{-1}(\gamma_j^{(i)})
     -	  	k_{ v(k)}  \eta_{i,v(k)}^{-1}(b) \Big)   + \nonumber \\
& +  E\Big(-\sum_{j=1}^{v(k)-1} k_j \eta_{i,j}^{-1}(\gamma_j^{(i)}) \Big)
   \sum_{r= v(k)+1}^{\tau_m} \sum_{b=0}^{\gamma_r^{(i)} -1} b^{-r}
      \nonumber \\
&   =  b^{-v(k)}   E\Big(-\sum_{j=1}^{v(k)-1} k_j \eta_{i,j}^{-1}(\gamma_j^{(i)}) \Big)
     \Big(   \sum_{b=0}^{\gamma_{ v(k)}^{(i)} -1}
     E \big(-k_{ v(k)} \eta_{i,v(k)}^{-1}(b)  \big)  \nonumber\\
& +  E\big( -k_{ v(k)} ( \eta_{i,v(k)}^{-1}(\gamma_{ v(k)}^{(i)})) \big)
   \{ b^{v(k)} [\gamma]_{\tau_m}^{(i)} \}
\Big) =    \hat{\d1}^{(i)}(k)   \nonumber.
\end{align}
Hence
\begin{equation} \nonumber
 \d1_{[0, [\gamma^{(i)}]_{\tau_m})}( x^{(i)}_{n})=\sum_{k_1,...,k_{\tau_m} \in \FF_b}  E(k \cdot   u^{(i)}_{\tilde{n}}) \hat{\d1}^{(i)}(k).
\end{equation}
Similarly, we obtain from \eqref{Lemm12a} and (\ref{Lemm12c}) that
\begin{equation} \nonumber
 \d1_{[0, \gamma^{(s+1)})}( x^{(s+1)}_{n})=\sum_{k_1,...,k_{m} \in \FF_b}  E(k \cdot   u^{(s+1)}_{\tilde{n}}) \hat{\d1}^{(s+1)}(k).
\end{equation}
Using (\ref{Lemm12}), we obtain
\begin{equation} \label{Lemm20}
  \prod_{i=1}^{s+1}  \d1_{[0, [\gamma^{(i)}]_{\tau_m})}( x^{(i)}_{n})  =
	\sum_{(k^{(1)},...,k^{(s)})  \in \FF_b^{\tau_m}}  \sum_{k^{(s+1)}  \in \FF_b^{m}}
	E(\mk \cdot   {\bf u}_{\tilde{n}}) 	\hat{\d1}(\mk).	
\end{equation}
Bearing in mind that $x^{(s+1)}_{n} =\{ n/b^m \}$ and $\gamma^{(s+1)}=N /b^m $, we have
\begin{equation} \nonumber
  \bx_{n +b^m A}  \in [\bs, [\bgamma]_{\tau_m}), n \in [0 ,N) \Longleftrightarrow
	 (\bx_{n +b^m A} , x^{(s+1)}_{n+b^m A}) \in [\bs, [\bgamma]_{\tau_m})  \times [0,\gamma^{(s+1)}).
\end{equation}
From \eqref{Lemm20} and \eqref{In1}, we derive
\begin{equation} \nonumber
     \Delta([\bs,[\bgamma]_{\tau_m}),(\bx_{n})_{n=b^m A}^{b^m A+N-1}) =
  \sum_{n=0}^{b^m-1} \prod_{i=1}^{s+1}   \d1_{[\bs, [\gamma^{(i)}]_{\tau_m})}( x^{(i)}_{n+b^m A})
	-b^m \prod_{i=1}^{s+1} [\gamma^{(i)}]_{\tau_m}
\end{equation}
\begin{equation} \nonumber
 =
  \sum_{n=0}^{b^m-1}	\;	\sum_{(k^{(1)}...,k^{(s)})  \in \FF_b^{\tau_m}} \;
   \sum_{k^{(s+1)}  \in \FF_b^{m}}
	 E(\mk \cdot \msu_{\widetilde{n+b^m A}}) \prod_{i=1}^{s+1} \hat{\d1}^{(i)}(k^{(i)})
		-b^m \prod_{i=1}^{s+1} [\gamma^{(i)}]_{\tau_m}  .
\end{equation}
Hence Lemma 2 is proved.   \qed \\


Let
\begin{equation} \nonumber
\mk = (k^{(1)},...,k^{(s+1)}),\;
 k^{(i)}=(k^{(i)}_{1},...,k^{(i)}_{\tau_m}), i \in [1,s], \;
k^{(s+1)}=(k^{(s+1)}_{1},...,k^{(s+1)}_{m}),
\end{equation}
\begin{equation} \nonumber
  G_m  =\{  \mk \; | \;  \;k^{(i)}_{j} \in \FF_b \; \with \;
  j \in [1,  \tau_m],\; i \in [1,s],\;\ad\;  j \in [1,  m] \;\for\; i=s+1 \},
\end{equation}	\\	
$G_m^{*} =  G_m \setminus \{ \bs\}$, and let
\begin{equation} \label{Lemm25}
     D_m =\{  \mk \in G_m \; | \;   \mk \cdot \msu_n  =0 \;\; \forall \;\;
						n \in [0,b^m-1] \}, \quad  \;\;  D_m^{*}  =  D_m \setminus \{ \bs\}.
\end{equation}						
It is easy to see that
\begin{equation} \label{Lemm25a}
   \mu \mk \in  D_m^{*} \quad \fall \quad \mu \in \FF_b^{*}, \; \mk \in  D_m^{*}.
\end{equation}	\\
{\bf Lemma 3.} {\it Let $(\bx_n)_{n \geq 0}$ be a digital sequence in base $b$. Then}
\begin{equation} \label{Lemm26}
     \Delta([\bs,[\bgamma]_{\tau_m}),(\bx_{n})_{n=b^m A}^{b^m A+N-1})
 =
    		\sum_{\mk \in G_m^{*}}   \hat{\d1}(\mk)  \sum_{n=0}^{b^m-1}
	 E(\mk \cdot \msu_n +\mk \cdot \msu_{ \widetilde{b^m A}}).
\end{equation} \\
{\bf Proof.} By  (\ref{Lemm12a}) we have $\hat{\d1}(\bs) =\prod_{i=1}^{s+1} [\gamma^{(i)}]_{\tau_m}$.
Applying Lemma 2, we get
\begin{equation} \nonumber
     \Delta([\bs,[\bgamma]_{\tau_m}),(\bx_{n})_{n=b^m A}^{b^m A+N-1})
 =
    		\sum_{\mk \in G_m^{*}}   \hat{\d1}(\mk)  \sum_{n=0}^{b^m-1}
	 E(\mk \cdot \msu_{ \widetilde{n+b^m A}}).
\end{equation}
 Using (\ref{Lemm2}), (\ref{Lemm3}) and (\ref{Lemm6}), we obtain
\begin{equation} \nonumber
    \widetilde{n+b^m A} =\widetilde{n} + \widetilde{b^m A} =\widetilde{n} \oplus \widetilde{b^m A}\quad \ad \quad
		\msu_{ \widetilde{n+b^m A}} =\msu_{ \widetilde{n}}
		+\msu_{ \widetilde{b^m A}}
\end{equation}
Now from  (\ref{Lemm26a}), we get
(\ref{Lemm26}). Hence Lemma 3 is proved. \qed \\ \\
{\bf Lemma 4.} {\it Let $(\bx_n)_{n \geq 0}$ be a digital sequence in base $b$. Then }
\begin{equation}   \label{In195}
  \sigma:= \sum_{n=0}^{b^m-1} E( \mk \cdot \msu_n ) =b^m \delta(\mk \in D_m).
\end{equation} \\
{\bf Proof.} Using  (\ref{Lemm3}) and (\ref{Lemm5}) and (\ref{Lemm12}) , we have
\begin{equation} \nonumber
 \mk \cdot \msu_{\tilde{n}} = \sum_{i=1}^s \sum_{j=1}^{\tau_m} \sum_{r=0}^{m-1}
	k_{j}^{(i)} \psi_r(a_r(n)) c_{j,r}^{(i)} + \sum_{j=1}^{m}
	 k_{j}^{(s+1)} \psi_{m-j}(a_{m-j}(n))
\end{equation}
\begin{equation}  \nonumber
 =\sum_{r=0}^{m-1} \psi_r(a_r(n)) \Big(\sum_{i=1}^s \sum_{j=1}^{\tau_m}
	k_{j}^{(i)} c_{j,r}^{(i)}
	 + k_{m-r }^{(s+1)} \Big) =\sum_{r=0}^{m-1} f_r \xi_r,
\end{equation}
where
\begin{equation}   \label{Lem4-1}
f_r=\psi_r(a_r(n)) \in \FF_b \quad \ad \quad
\xi_r=\sum_{i=1}^s \sum_{j=1}^{\tau_m}
	k_{j}^{(i)} c_{j,r}^{(i)}
	 + k_{m-r }^{(s+1)}.
\end{equation}
By (\ref{Lemm26a}), (\ref{In16}) and (\ref{Del2}), we obtain
\begin{equation}   \nonumber
  \sigma = \sum_{\tilde{n}=0}^{b^m-1} E( \mk \cdot \msu_{\tilde{n}}) =\sum_{f_0,...,f_{m-1} \in \FF_b} E(\sum_{r=0}^{m-1} f_r \xi_r)=
	b^m \prod_{r=0}^{m-1} \delta (\xi_r =0).
\end{equation}
Now from (\ref{Lemm25}), we get that $ \mk \in D_m$ and Lemma 4 follows. \qed  \\

Let
\begin{equation}   \nonumber
  \Lambda_m=\{\mk =(k^{(1)},...,k^{(s+1)}) \in G_m \; | \; k^{(s+1)}=\bs\} ,
\end{equation}
\begin{equation}   \nonumber
  g_{\bw} = \{ A \geq 1 \;  | \;  y^{(i)}_{b^m A,j}= w^{(i)}_j, \;i \in [1,s],\;j \in
	[1,\tau_m]  \},   \quad  \rho_{\bw}= 0 \;\; \for \;\; g_{\bw} = \emptyset,
\end{equation}
\begin{equation}   \label{End0}
 \rho_{\bw} = \min_{A \in g_{\bw} }   A \;\; \for \;\; g_{\bw} \neq \emptyset,
\quad \;\;
   M_m =\{\rho_{\bw} \; | \; \bw \in \Lambda_m  \}.
\end{equation}
We consider the following conditions :
\begin{equation}    \label{End}
 g_{\bw} \neq \emptyset \quad \fall \quad \bw \in \Lambda_m
\end{equation}
and
\begin{equation}    \label{End1}
\sigma_1:=  \frac{1}{\card(R_m)} \sum_{A \in R_m}
| \Delta([\bs,[\bgamma]_{\tau_m}),(\bx_{n})_{n=b^m A}^{b^m A+N-1})|^2
 =\sum_{\mk  \in  D_m^{*}}  b^{2m}  |\hat{\d1} (\mk  )|^2
\end{equation}
for some finite set $R_m$.\\
Bearing in mind (\ref{Lemm3}), we get
\begin{equation}   \label{End2}
  g_{\bw} = \{ A \geq 1 \;  | \;  u^{(i)}_{\widetilde{b^m A }}=w^{(i)}, \; i \in [1,s]\},
	\;\; \where \;\; w^{(i)}=(w^{(i)}_1,...,w^{(i)}_{\tau_m}).
\end{equation}\\
{\bf Lemma 5.} {\it   Let $(\bx_n)_{n \geq 0}$ be a weakly admissible uniformly distributed digital
 $(\bT,s)$-sequence in base $b$, satisfying to (\ref{End}) for  all $m \geq m_0$ with some $m_0 \geq 1$.
 Then (\ref{End1}) is true for $R_m=M_m$}.\\ \\
{\bf Proof.}
By  (\ref{End0}) and (\ref{Del2}), we obtain
\begin{equation}   \nonumber
\frac{1}{b^{s \tau_m}} 	\sum_{\bw \in \Lambda_m}  E(\mk  \cdot \bw) =
\frac{1}{b^{s \tau_m}} 	\sum_{w^{(i)}_j \in \FF_b, \;i\in[1,s],\; j \in [1,\tau_m] }
 E\Big(\sum_{i=1}^s \sum_{j=1}^{\tau_m}
k^{(i)}_j w^{(i)}_j \Big)
\end{equation}
\begin{equation}   \label{In200}
= \prod_{i=1}^s\prod_{j=1}^{\tau_m} \delta(k^{(i)}_j=0)
= \prod_{i=1}^s \delta(k^{(i)}=0) ,\quad \where \quad \mk \in G_m .
\end{equation}
Using (\ref{Lemm26}), we derive
\begin{equation} \label{In201}
 | \Delta([\bs,[\bgamma]_{\tau_m}),(\bx_{n})_{n=b^m A}^{b^m A+N-1})|^2
\end{equation}	
\begin{equation} \nonumber
 =	      \sum_{\dot{\mk}, \ddot{\mk} \in G_m^{*}} \hat{\d1}(\dot{\mk} )
				  \overline{\hat{\d1}(\ddot{\mk} ) }
				\sum_{\dot{n},\ddot{n}=0}^{b^{m}-1}
E( \dot{\mk} \cdot \msu_{\dot{n} } +
   \dot{\mk} \cdot \msu_{\widetilde{b^m A}}
             -\ddot{\mk} \cdot \msu_{\ddot{n} }
						-\ddot{\mk} \cdot \msu_{\widetilde{b^m A}}).
\end{equation}
It is easy to see that if condition (\ref{End}) is true, than $\card( M_m) =b^{s \tau_m}$ and
 $\{ \msu_{\widetilde{b^m A}} \; | \; A \in M_m \} =\Lambda_m$. \\
Applying  (\ref{End0}), (\ref{In200}), (\ref{In201}) and (\ref{End1}) with $R_m =M_m$, we have
\begin{equation} \nonumber
\sigma_1=
	      \sum_{\dot{\mk}, \ddot{\mk} \in G_m^{*}} \hat{\d1}(\dot{\mk} )
				  \overline{\hat{\d1}(\ddot{\mk} )}
				\sum_{\dot{n},\ddot{n}=0}^{b^{m}-1}
				b^{-s\tau_m}  \sum_{A \in M_m}
E( \dot{\mk} \cdot \msu_{\dot{n}} -
             \ddot{\mk} \cdot \msu_{\ddot{n}}
						+(\dot{\mk} -\ddot{\mk}) \cdot \msu_{\widetilde{b^m A}} )
\end{equation}
\begin{equation} \nonumber
 =
	      \sum_{\dot{\mk},  \ddot{\mk} \in G_m^{*}} \hat{\d1}(\dot{\mk} )
				  \overline{\hat{\d1}(\ddot{\mk} )}
				\sum_{\dot{n},\ddot{n}=0}^{b^{m}-1} 	
E( \dot{\mk} \cdot \msu_{\dot{n}} -
             \ddot{\mk} \cdot \msu_{\ddot{n}}  )
				b^{-s\tau_m}		\sum_{\bw \in\Lambda_m} E(\dot{\mk} -\ddot{\mk}) \cdot \bw  )
\end{equation}
\begin{equation} \nonumber
 =
	     \sum_{\dot{\mk}, \ddot{\mk} \in G_m^{*}} \hat{\d1}(\dot{\mk} )
				  \overline{\hat{\d1}(\ddot{\mk} ) }
				\sum_{\dot{n},\ddot{n}=0}^{b^{m}-1} 	
E( \dot{\mk} \cdot \msu_{\dot{n}} -
             \ddot{\mk} \cdot \msu_{\ddot{n}}  )
				\prod_{i=1}^s \delta( \dot{k}^{(i)} = \ddot{k}^{(i)}).
\end{equation}
%
%
Let $\dddot{n} =\ddot{n} \ominus \dot{n}$.
From (\ref{In16}), we obtain $\{ \dddot{n}\;| 0 \leq \ddot{n} <b^m\} = \{0,1,...,b^m-1 \}$.\\
By (\ref{Lemm2}) - (\ref{Lemm6}), we get $\msu_{\dddot{n}} =\msu_{\ddot{n}}-\msu_{\dot{n}}$.
Hence
\begin{equation} \nonumber
 \sigma_1=
	      \sum_{\dot{\mk}, \ddot{\mk} \in G_m^{*}} \hat{\d1}(\dot{\mk} )
				  \overline{\hat{\d1}(\ddot{\mk} )}
				\sum_{\dot{n},\dddot{n}=0}^{b^{m}-1}
E( (\dot{\mk} - \ddot{\mk})\cdot \msu_{\dot{n}} -
\ddot{\mk} \cdot \msu_{\dddot{n}})       \prod_{i=1}^s \delta( \dot{k}^{(i)} - \ddot{k}^{(i)}).
\end{equation}
We get
 $\dot{\mk} - \ddot{\mk}=(0,...,0,\dot{k}^{(s+1)}-\ddot{k}^{(s+1)})$.\\
From (\ref{Lemm5}), we have
$u^{(s+1)}_{\tilde{n},j} =\omega^{-1}(a_{m-j+1}(n))$ and
\begin{equation}   \nonumber
  (\dot{\mk} - \ddot{\mk}) \cdot \msu_{\tilde{n}}=
	(\dot{k}^{(s+1)}-\ddot{k}^{(s+1)}) \cdot u^{(s+1)}_{\tilde{n}}
= \sum_{j=1}^m  (\dot{k}^{(s+1)}_{j} -\ddot{k}^{(s+1)}_{j}) \omega^{-1}(a_{m-j}(\tilde{n})).
\end{equation}
Taking into account (\ref{Del2}),  we get
\begin{equation}   \nonumber
\sum_{\dot{n}=0}^{b^{m}-1}  E(
  (\dot{\mk} - \ddot{\mk}) \cdot \msu_{\dot{n}})
 	   = \sum_{\tilde{n}=0}^{b^{m}-1}  E(
  (\dot{\mk} - \ddot{\mk}) \cdot \msu_{\tilde{n}}) =
b^m \prod_{j=1}^m \delta(\dot{k}^{(s+1)}_j=\ddot{k}^{(s+1)}_j) .
\end{equation}
Hence $\dot{\mk}=  \ddot{\bk}$. Using Lemma 4, we obtain
\begin{equation} \nonumber
 \sigma_1=
b^m\sum_{\dot{\mk} \in  G_m^{*} }|\hat{\d1}(\dot{\mk} ) |^2
				\sum_{\dddot{n}=0}^{b^{m}-1}
E(-\dot{\mk} \cdot \msu_{\dddot{n}} ) =\sum_{\mk  \in  D_m^{*}}  b^{2m}  |\hat{\d1} (\mk  )|^2.
\end{equation}
Therefore Lemma 5 is proved. \qed  \\ \\
Let  $\Psi$ be a set of all bijections $ \dot{\psi} : \ZZ_b \to  \FF_b $,
 $\psi \in \Psi$,
$k \in \FF_b,\\ c \in \ZZ_b $,
   $A_{k,c, \psi}=
E(-k \psi(c) )\sum_{b=0}^{c-1} E(k \psi(b) )$, $\;\langle x	\rangle =\min(\{x\},1- \{x\})$,  $x \in [0,1]$ and let
\begin{equation}   \label{Lem2a}
B_{k,c, \psi}(x)= \sum_{b=0}^{c-1} E(k \psi(b) )  +E(k \psi(c))x
=E(k \psi(c)  /q)(A_{k,c, \psi} +x) .
\end{equation}
 By (\ref{Lemm12a}) and (\ref{Lem2a}), we get
\begin{equation}   \label{Lem2d}
 b^{v(\dot{k})} |\hat{\d1}^{(s+1)}  (\mu \dot{k}) | = |B_{\tilde{k},\tilde{c},\psi_1}(x_1)| \qquad
  \ad \qquad   b^{v(\ddot{k})} |\hat{\d1}^{(s)}  (\mu \ddot{k}) | = |B_{\breve{k}   ,\breve{c},\psi_2}(x_2)| .
\end{equation}
with $\tilde{k} =-\mu \dot{k}_{v(\dot{k})}$, $\breve{k} =-\mu \ddot{k}_{v(\ddot{k})}$,
 $\tilde{c} = \gamma^{(s+1)}_{v(\dot{k})}  $,   $\breve{c} = \gamma^{(s)}_{v(\ddot{k})}  $,
  $ \psi_1 = \psi^{-1}_{v(\dot{k})} $,
  $ \psi_2 = \eta^{-1}_{s,v(\ddot{k})} $, \\  $x_1 =\{ b^{v(\dot{k})} \gamma^{(s+1)} \}$,
 $x_2 =\{ b^{v(\ddot{k})} [\gamma^{(s)}]_{\tau_m} \}$.  \\ \\
{\bf Lemma 6.} {\it  With the notations as above, there exist $\sa_1,...,\sa_{b+7}\in \ZZ_b$, $\sa_1^2+...+\sa_{b+5}^2 >0$, $\sa_{b+6} = \sa_{b+7} =0$ such that}
\begin{equation}   \label{Lem2}
           \Big| B_{k,c, \psi}\Big(\sum_{j=1}^{b+7}\frac{\sa_j}{b^j}+\frac{y}{b^{b+7}}\Big) \Big| \geq b^{-b-7}, \;\;\;
			\forall \;			k \in \FF_b,  c \in \ZZ_b,  y\in [0,1], \psi \in \Psi
\end{equation}
and
\begin{equation}   \label{Lem3}
  \sum_{k \in \FF_b^{*}}  |B_{k,c, \psi}(x)|^2  \geq b^{-2r}  \;\; \forall \;		  c \in \ZZ_b,	   \;\; \where
   \;\;  \langle x	\rangle \geq b^{-r}, \;   r \geq 1 .
\end{equation} \\
{\bf Proof.} Let
\begin{equation}   \nonumber
  \dot{A}= \{  \theta_{k,c, \psi}:=\Re(A_{k,c, \psi}) \;|\; k \in \FF_b,\;  c \in \ZZ_b, \; \psi \in \Psi\}.
\end{equation}
Taking into account that $\card(\Psi) =b!$, we get $\card(\dot{A}) \leq b!b^2+2< b^{b+4}$.\\
Let
\begin{equation}   \nonumber
 \ddot{A}= \{ \ba =(\sa_1,...,\sa_{b+7}) \in \ZZ_b^{b+7}  \;|\;
\sa_1^2+...+\sa_{b+5}^2 >0, \quad \sa_{b+6} = \sa_{b+7} =0 \}
\end{equation}
and let $z_{\ba}=\sa_1/b+ \cdots + \sa_{b+7}/b^{b+7} $.
By	 (\ref{Lem2a}), we derive
\begin{equation}   \nonumber
          |B_{k,c,\psi}(x)| =|A_{k,c,\psi} +x| \geq |\Re(A_{k,c,\psi}) +x|.
\end{equation}
Suppose that  (\ref{Lem2}) is not true. Then for all $\ba \in \ddot{A}$ there exist $k(\ba),c(\ba),\psi(\ba)$ and $y(\ba)$
such that
\begin{equation}   \nonumber
  b^{-b-7} >   \Big|B_{k(\ba),c(\ba),\psi(\ba)}\Big(\sum_{j=1}^{b+7}\frac{\sa_j}{b^j}+\frac{y(\ba)}{b^{b+7}}\Big) \Big| \geq
  \Big|\theta_{k(\ba),c(\ba),\psi(\ba)} +z_{\ba} +\frac{y(\ba)}{b^{b+7}} \Big|.
\end{equation}
Hence $ |\theta_{k(\ba),c(\ba),\psi(\ba)} +z_{\ba}  | <   b^{-b-6}  $.
Suppose that $\theta_{k(\ba_1),c(\ba_1),\psi(\ba_1)}  = \theta_{k(\ba_2),c(\ba_2),\psi(\ba_2)} $ for some $\ba_1 , \ba_2 \in \ddot{A} $, $\ba_1  \neq \ba_2 $. Hence $|z_{\ba_1} -z_{\ba_2}| < b^{-b-5} $.
Bearing in mind that $|z_{\ba_1} -z_{\ba_2}| \geq b^{-b-5} $ for all $\ba_1 \neq \ba_2$, we
get a contradiction. Therefore  $\theta_{k(\ba_1),c(\ba_1),\psi(\ba_1)}  \neq \theta_{k(\ba_2),c(\ba_2),\psi(\ba_2)} $ for all $\ba_1 , \ba_2 \in \ddot{A} $, $\ba_1  \neq \ba_2 $.
Thus $ \card(\dot{A}) \geq \card(\ddot{A})$.
Hence
\begin{equation}   \nonumber
        b^{b+4} > \card(\dot{A}) \geq \card(\ddot{A}) =b^{b+5}-1 >  b^{b+4}.
\end{equation}
We have a contradiction. Therefore (\ref{Lem2}) is true.

Now we consider assertion (\ref{Lem3}).
If $c=0$, then  $|B_{k,c, \psi}(x)|=|x|$ and (\ref{Lem3}) follows. \\
 Now let $c \in \{ 1,...,b-1\}$.
 By  (\ref{Lem2a}), we have
\begin{equation}   \label{Lem4}
 |B_{k,c, \psi}(x)|^2 =|A_{k,c, \psi}|^2 + x (A_{k,c, \psi} + \overline{A_{k,c, \psi}}) +x^2.
\end{equation}
Using (\ref{Del2}), we get
\begin{equation}   \nonumber
  \sum_{k \in \FF_b^{*}}  |A_{k,c, \psi}|^2 =-c^2 + \sum_{b_1,b_2=0}^{c-1} 	
	\sum_{k \in \FF_b } E(k (\psi(b_1)-\psi(b_2))
\end{equation}
\begin{equation}   \nonumber
 =-c^2 +b\sum_{b_1,b_2=0}^{c-1}
	\delta(\psi(b_1)=\psi(b_2))  = -c^2 +bc.
\end{equation}
Taking into account that $ \psi(0) =0$, we obtain
\begin{equation}   \nonumber
  \sum_{k \in \FF_b^{*}}  A_{k,c, \psi} =-c + \sum_{b=0}^{c-1}
	\sum_{k \in \FF_b } E(k \psi(b) )=-c +b \sum_{b=0}^{c-1}  \delta( \psi(b) =0 ) \geq b-c.	
\end{equation}
Now from (\ref{Lem4}), we derive
\begin{equation}   \nonumber
  \sum_{k \in \FF_b^{*}}  |B_{k,c, \psi}(x)|^2  \geq  c(b-c) +x^2 +2x(b-c) \geq x^2
\end{equation}
and (\ref{Lem3}) follows. Thus Lemma 6 is proved.    \qed  \\ \\
Applying (\ref{Lem2d}) - (\ref{Lem3}),   we have  \\ \\
{\bf Corollary.} {\it Let  $\sa_1,...,\sa_{b+7}$ be integers chosen in Lemma 6 and let
$\gamma^{(s+1)}_{v(\dot{k})+ j} =\sa_j$,   $j=1,...,b+7$, with some $\dot{k} $.
%
Then}
\begin{equation}   \label{Lem2ab}
 		|\hat{\d1}^{(s+1)}  (\mu \dot{k}) | \geq b^{-v(\dot{k})-b-7}  \;\;\;
			\forall \;			\mu \in \FF_b^{*}
\end{equation}
and
\begin{equation}   \label{Lem3b}
   \sum_{\mu \in \FF_b^{*} } |\hat{\d1}^{(i)}  (\mu \ddot{k}) |^2 \geq b^{-2v(\ddot{k})-2\dot{r}} , \;\;
	\where \;\; \; \langle b^{v(\ddot{k})} \gamma^{(i)}	\rangle \geq b^{-\dot{r}}.
\end{equation}\\
{\bf Lemma 7.} {\it Let $(\bx_n)_{n \geq 0}$ be a digital sequence in base $b$  and
 let  $\rho \in[2,m-2]$ be an integer.
 Then there exists $ \mk \in D_m^{*}$ such that $k^{(1)}=... =k^{(s-1)}=0$, $ k^{(s)}_{ v(k^{(s)})}=1$,
$1 \leq v(k^{(s)}) \leq \rho-1 $ and $v(k^{(s+1)}) \leq m-\rho+2$ }. \\ \\
{\bf Proof.} From (\ref{Lemm2})-(\ref{Lemm5}), (\ref{Lemm25}) and (\ref{Lem4-1}), we get that $ \mk \in D_m^{*}$ if and only if
\begin{equation}   \label{Lem61}
 \sum_{i=1}^{s}\sum_{j=1}^{\tau_m}
	k_{j}^{(i)} c_{j,r}^{(i)}
	 + k_{m-r }^{(s+1)} = 0, \quad \for \; {\rm all} \quad r=0,1,...,m-1.
\end{equation}
We put $k^{(1)}=... =k^{(s-1)}=0, \; k^{(s)}_{j} = 0$, for $j \geq \rho$ and $ k^{(s+1)}_{j} = 0$,
 for $j > m-\rho +2$.
Hence (\ref{Lem61}) is true  if and only if
\begin{equation}   \label{Lem6a}
 k_{m-r }^{(s+1)}  = -\sum_{j=1}^{\rho-1}
	k_{j}^{(s)} c_{j,r}^{(s)}
	\;\; \for  \;\; r=0,1,...,m -1, \quad k_{m-r }^{(s+1)}  = 0 \; \for \; m-r >  m-\rho +2.
\end{equation}
Therefore, in order to obtain the statement of the lemma, it is sufficient to show that there exists a nontrivial solution of the following system of linear equations
\begin{equation}   \label{Lem6}
 \sum_{j=1}^{\rho-1}
	k_{j}^{(s)} c_{j,r}^{(i)}
	 + k_{m-r }^{(s+1)} \delta(m-r \leq m-\rho+2)=0,   \quad  r=0,...,m-1.
\end{equation}
 In this system, we have  $m+1$ unknowns
$k_{1}^{(s)},...,k_{\rho-1}^{(s)}$, $k_{1}^{(s+1)},...,k_{m-\rho+2}^{(s+1)}$
and $m$ linear equations.  Hence there exists a nontrivial solution of (\ref{Lem6}).
 By  (\ref{Lem6}), we get that if $k^{(s)} =0$, then $k^{(s+1)} =0$. Hence $k^{(s)} \neq 0$ and
$1 \leq v(k^{(s)}) \leq \rho-1 $.
Taking into account  that if $\mk \in D_m$ then $ \mu \mk \in D_m$ for all $\mu \in \FF_b^{*}$.
Therefore there exists $ \mk \in D_m^{*}$ such that  $ k^{(s)}_{ v(k^{(s)})}=1$ and
$1 \leq v(k^{(s)}) \leq \rho-1 $.
 Thus Lemma 7 is proved. \qed \\ \\
{\bf Proposition.}  {\it   Let $(\bx_n)_{n \geq 0}$ be a weakly admissible uniformly distributed
  digital $(\bT,s)$-sequence in base $b$, satisfying to (\ref{End}) for  all $m \geq m_0 \geq 1$.
Then  $[0,\gamma_1) \times ...\times [0,\gamma_s)$ is of bounded remainder  with respect to
  $(\bx_n)_{n \geq 0}$ if and only if (\ref{Cond}) is true.}\\  \\
{\bf Proof.}  The sufficient part  of the Theorem and of the Proposition follows directly from the definition of $(\bT, s)$
  sequence and Lemma B. We will consider only the necessary part of the Theorem and of the Proposition.
	
	Suppose that  (\ref{Cond}) does not true. Then
\begin{equation} \nonumber
   \max_{1 \leq i \leq s} \card \{ j \geq 1 \; | \;\gamma_{j}^{(i)} \neq 0 \}  = \infty.
\end{equation}
Let, e.g.,
\begin{equation} \nonumber
   \card \{ j \geq 1 \; | \; \gamma_{j}^{(s)} \neq 0 \}  = \infty.
\end{equation}
Let
\begin{equation} \label{Prop0}
 W = \{ j \geq 1 \; | \; \gamma_{j}^{(s)} \in \{ 1,...,b-2\}   \;\oor\;
\gamma_{j}^{(s)} =b-1, \; \and  \gamma_{j+1}^{(s)} =0 \}   .
\end{equation}
Bearing in mind that $ \{ j \geq 1 \; | \;\gamma_{l}^{(s)} = b-1 \;\forall \; l >j \}   =\emptyset$,
 we obtain  that $\card (W)=\infty$.

Suppose that there exists $H>0$ such that  $b^{2H}c_1 > 4H^2$, $ c_1 = \gamma_0^2 b^{-4b-36}  $,
\begin{equation} \label{Prop1}
   |\Delta([\bs,\bgamma),(\bx_{n})_{n=M}^{M+N-1})| \leq H -s
	\quad \fall \quad   \;M\geq 0, \; N \geq 1,
\end{equation}
with $[\bs,\bgamma) =[0,\gamma_1)\times \cdots \times [0,\gamma_s)$,
$ \gamma_0 = \gamma_1 \gamma_2 \cdots \gamma_{s-1} $.

Let $W=\{\dot{w}_j \; | \; \dot{w}_i <\dot{w}_j \; \for \; i<j,\; \; j=1,2,...\}$  and let
\begin{equation} \label{Prop2}
  r(1)=\dot{w}_1,\;\;  r(j+1) = \min(\dot{w}_k \in W \; | \; \dot{w}_k \geq r(j) +H^2 ), \;\; j=1,2,...\;.
\end{equation}
We choose $m$ and $J$ from the following conditions
\begin{equation} \label{Cond2}
  m=r(J)+b+10, \;  2\prod_{i=1}^{s-1} [\gamma_i]_{\tau_m} \geq
	    \prod_{i=1}^{s-1}\gamma_i =\gamma_0,  \;  J \geq H^2 b^{2b+30}\gamma_0^{-2} ,
  \; m \geq m_0.
\end{equation}
Applying Lemma 1 and (\ref{Prop1}), we have
\begin{equation} \label{Prop3}
   |\Delta([\bs,[\bgamma]_{\tau_m}),(\bx_{n})_{n=b^m A}^{b^m A+N-1})| \leq H
	\qquad  \forall \;A \geq 0, \; N \in [1, b^m].
\end{equation}
By Lemma 7, we get that there exists a sequence $(\mk(j))_{j=1}^J$ such that
\begin{equation} \nonumber
 \mk(j) \in D_m^{*}, \qquad k^{(1)}(j)=...=k^{(s-1)}(j)=0, \;\; k^{(s)}_{ v(k^{(s)}(j))}(j)=1,
\end{equation}
\begin{equation} \label{Prop4}
  v(k^{(s)}(j)) \leq r(j)-1, \qquad v(k^{(s+1)}(j)) \leq m-r(j)+2, \;\;\; j \in [1,J].
\end{equation}
We see that the  sequence $(\mk(j))_{j=1}^J$ does not depend on $\gamma^{(s+1)}$.\\
Using (\ref{Prop0}) and (\ref{Prop2}), we obtain $\gamma^{(s)}_{r(j)} \neq 0$. Hence
\begin{equation}   \label{Prop4a}
   \langle b^{{v(k^{(s)}(j))}}\gamma^{(s)}\rangle 	=.\gamma^{(s)}_{v(k^{(s)}(j))+1}
    \cdots \gamma^{(s)}_{r(j)} \cdots   \; \geq \; b^{v(k^{(s)}(j))-r(j)-2},
\end{equation}
$j=1,...,J$.
Let $H_1 = \{  1,2,...,J\}$ if
\begin{equation} \label{Prop5}
       |v(k^{(s+1)}(j)) - v(k^{(s+1)}(j_1))| \geq b+8
\end{equation}
for all $1 \leq j < j_1 \leq J$, and let
 $H_1 = \{  j\}$ if there exist $1 \leq j < j_1 \leq J$ such that \eqref{Prop5} is false.
 Let  $\sa_1,...,\sa_{b+7}$ be integers chosen in Lemma 6  and let
\begin{equation} \label{Prop8}
 N =b^m \gamma^{(s+1)}  \quad \with \quad \gamma^{(s+1)} =\sum_{j \in H_1}
\sum_{\nu=1}^{b+7}\sa_{\nu} b^{\nu+v(k^{(s+1)}(j))}.
\end{equation}
%
From Lemma 5, (\ref{End0}), (\ref{End1}),  (\ref{Prop3}) and conditions of the Proposition, we have
\begin{equation}   \label{Prop9}
H^2 \geq \sigma_1 = \sum_{\mk \in D_m^{*}}   b^{2m} |\hat{\d1} (\mk)|^2.
\end{equation}
Taking into account (\ref{Lemm25a}), (\ref{Prop4}) and (\ref{Prop5}), we get that
if $\mk(j) \in D_m$ then $ \mu \mk(j) \in D_m$ for $\mu \in \FF_b^{*}$,
and if  $j_1,j_2 \in H_1$, $j_1 \neq j_2$, then     $ \mu_1 \mk(j_1) \neq \mu_2 \mk(j_2) $ for all
 $\mu_1,\mu_2 \in \FF_b^{*}$.\\
According to (\ref{Lemm12}), (\ref{Lemm12a}) and (\ref{Prop9}), we have
\begin{equation}   \nonumber
         \sigma_1 \geq
   \sum_{\mu \in \FF_b^{*} } \sum_{j \in H_1}  b^{2m}  |\hat{\d1} (\mu \mk(j))|^2
\end{equation}
\begin{equation}   \nonumber
              =([\gamma_1]_{\tau_m} \cdots [\gamma_{s-1}]_{\tau_m})^2
   \sum_{\mu \in \FF_b^{*} } \sum_{j \in H_1}  b^{2m}  |\hat{\d1}^{(s)} (\mu k^{(s)}(j))|^2
  |\hat{\d1}^{(s+1)} (\mu k^{(s+1)}(j))|^2  .
\end{equation}
From Corollary  and (\ref{Prop8}), we obtain
\begin{equation}   \nonumber 
      |\hat{\d1}^{(s+1)} (\mu k^{(s+1)}(j))|^2 \geq b^{-2v(k^{(s+1)}(j))-2b-14} \qquad \fall\;  \mu \in \FF_b^{*},\; j \in H_1.
\end{equation}
By (\ref{Prop4a}), we can apply   Corollary with $ \dot{r} =r(j) - v(k^{(s)}(j)) +2$. Hence
\begin{equation}   \label{Lem2c}
         \sum_{\mu \in \FF_b^{*} } |\hat{\d1}^{(s)}(\mu k^{(s)}(j)) |^2
				\geq b^{-2v(k^{(s)}(j))  -2( r(j)-v(k^{(s)}(j)) +2)}    = b^{-2r(j)-4}, \quad j \in H_1.
\end{equation}
Using (\ref{Prop9})-(\ref{Lem2c}) and (\ref{Cond2}), we obtain
\begin{equation}   \nonumber
  4H^2 \geq 4\sigma_1 \geq  \sigma_1 \gamma_0^{2} ([\gamma_1]_{\tau_m} \cdots [\gamma_{s-1}]_{\tau_m})^{-2}
\geq  \gamma_0^{2} \sum_{j \in H_1}
 \sum_{\mu \in \FF_b^{*} }
 |\hat{\d1}^{(s+1)} (\mu k^{(s)}(j))|^2
\end{equation}
\begin{equation}   \label{In97}
  \times b^{2m -2v(k^{(s+1)}(j))-2b-14}   \geq \gamma_0^{2} \sum_{j \in H_1}
 b^{2m- 2r(j) -2v(k^{(s+1)}(j))-2b-18 }.
\end{equation}
Suppose that $\card(H_1) =J$.
From (\ref{Cond2}) and (\ref{Prop4}), we get
\begin{equation}   \label{In98}
 4H^2 \geq 4\sigma_1 \geq \gamma_0^{2} \sum_{j=1}^{J}
  b^{2m- 2r(j) -2v(k^{(s+1)}(j))-2b-18 } \geq \gamma_0^{2} J b^{-2b-22} >4H^2.
\end{equation}
We have a contradiction. Now let $\card(H_1) =1$.

By (\ref{Prop5}), we obtain that there exist $j,j_1 \in [1,J]$ such that $j \in H_1$, $j <j_1$ and
 $|v(k^{(s+1)}(j))-v(k^{(s+1)}(j_1))| \leq b+7$.

According to  (\ref{Prop2}) and (\ref{Prop4}), we have
\begin{align}
& r(j) +v(k^{(s+1)}(j)) \leq   r(j_1) -H^2 + v(k^{(s+1)}(j_1))+b+7 \;\; \ad \;m - r(j) \nonumber \\
&   - v(k^{(s+1)}(j)) \geq  m - r(j_1)  - v(k^{(s+1)}(j_1))+H^2 -b -7  \geq H^2 -b-9.  \nonumber
\end{align}
Applying  \eqref{Prop1}, \eqref{Prop3} and \eqref{In97}, we get
\begin{equation}   \nonumber
 4H^2 \geq 4 \sigma_1 \geq \gamma_0^{2}  b^{2m- 2r(j) -2v(k^{(s+1)}(j))-2b-18 } \geq
 \gamma_0^{2}  b^{2H^2-4b-36}=b^{2H^2}c_1 >4H^2 ,
\end{equation}
with $ c_1 = \gamma_0^2 b^{-4b-36}  $. We have a contradiction.
By  \eqref{In98}, the Proposition is proved. \qed \\ \\
{\bf Completion of the proof of the Theorem.}
 By Lemma A,
 $S(\bL^{(m)})$ is a  uniformly distributed digital $(T,s)$-sequence in base $b$.

By Theorem A, we get that
 $1, L_1, . . . , L_s$ are
linearly independent over $\FF_b[z]$. Hence  $1, z^m L_1, ...,
z^m L_s$ are
linearly independent over $\FF_b[z]$. Let
 $\bL^{(m)} = (z^m L_1, . . . ,z^m L_s)$, and  let  $S(\bL^{(m)}) = (\bl^{(m)}_n)_{n \geq 0}$
 (see (\ref{In22}))
with
\begin{equation} \nonumber
\bl_n^{(m)} = (l^{(m,1)}_n,..., l^{(m,s)}_n   ), \quad
l^{(m,i)}_n = \eta^{(i)}( n(z)z^m L_i(z)),
 \quad \for \quad  \; 1 \leq i \leq s, \; n \geq 0.
\end{equation}
 Using Theorem A, we obtain that $S(\bL^{(m)})$ is a uniformly distributed sequence in $[0,1)^s$.
Therefore, for all $\bw \in \Lambda_m $ there exists an integer $A \geq 1$ with
\begin{equation} \nonumber
l_{b^m A,j}^{(m,i)}=\eta^{(i)}_j(w^{(i)}_j)
 \quad \for \quad  \; 1 \leq i \leq s, \; 1 \leq j \leq \tau_m.
\end{equation}
Thus  $S(\bL^{(m)})$ satisfies the condition (\ref{End}).

Bearing in mind that  $1, L_1, . . . , L_s$ are
linearly independent over $\FF_b[z]$, we get  that $\{ n(z) L_i\} \neq 0$ for all $n \geq 1$. Hence  $\{  l^{(i)} (n)\} \neq 0$ for all $n \geq 1$
$(i=1,...,s)$. Therefore the sequence   $S(\bL)$ is weakly admissible.

 Applying the Proposition, we get the assertion of the Theorem. \qed

{\bf Address}: Department of Mathematics,
Bar-Ilan University, Ramat-Gan, 5290002, Israel \\
{\bf E-mail}: mlevin@math.biu.ac.il\\

\end{document}